\documentclass{amsart}
\usepackage{amssymb,curves}
\newtheorem{Lemma}{Lemma}[section]\newcommand{\bel}{\begin{Lemma}}\newcommand{\eel}{\end{Lemma}}
\newtheorem{Proposition}[Lemma]{Proposition}\newcommand{\bprop}{\begin{Proposition}}\newcommand{\eprop}{\end{Proposition}}
\newtheorem{Theorem}[Lemma]{Theorem}\newcommand{\bthe}{\begin{Theorem}}\newcommand{\ethe}{\end{Theorem}}
\newcommand{\bpr}{{\it Proof:~}}\newcommand{\epr}{$\blacksquare$\\}
\newtheorem{Remark}[Lemma]{Remark}\newcommand{\beR}{\begin{Remark}\rm}\newcommand{\eeR}{\end{Remark}}
\newtheorem{Definition}[Lemma]{Definition}\newcommand{\bed}{\begin{Definition}}\newcommand{\eed}{\end{Definition}}
\newtheorem{Example}[Lemma]{Example}\newcommand{\bex}{\begin{Example}\rm}\newcommand{\eex}{\end{Example}}
\newtheorem{Corollary}[Lemma]{Corollary}\newcommand{\bcor}{\begin{Corollary}\rm}\newcommand{\ecor}{\end{Corollary}}
\newtheorem{Fact}[Lemma]{Fact}\newcommand{\bfact}{\begin{Fact}\rm}\newcommand{\efact}{\end{Fact}}

\newcommand{\beq}{\begin{equation}}\newcommand{\eeq}{\end{equation}}
\newcommand{\bem}{\begin{displaymath}}\newcommand{\eem}{\end{displaymath}}
\newcommand{\beqa}{\begin{eqnarray}}\newcommand{\eeqa}{\end{eqnarray}}
\newcommand{\bee}{\begin{enumerate}}\newcommand{\eee}{\end{enumerate}}
\newcommand{\bei}{\begin{itemize}}\newcommand{\eei}{\end{itemize}}
\newcommand{\bet}{\begin{tabular}{cccccccc}}\newcommand{\eet}{\end{tabular}}
\newcommand{\bpm}{\begin{pmatrix}}\newcommand{\epm}{\end{pmatrix}}
\newcommand{\bM}{\begin{matrix}}\newcommand{\eM}{\end{matrix}}
\newcommand{\ber}{\begin{array}{l}}\newcommand{\eer}{\end{array}}

\newcommand{\tinyM}{\scriptstyle}
\newcommand{\tinyT}{\scriptsize}\newcommand{\smallT}{\footnotesize}
\newcommand{\tinyA}{\tinyM\text{\tinyT}}

\newcommand{\mC}{\mathbb{C}}

\newcommand{\mR}{\mathbb{R}}

\newcommand{\mZ}{\mathbb{Z}}

\newcommand{\Ga}{\Gamma}
\newcommand{\ep}{\epsilon}\newcommand{\si}{\sigma}

\newcommand{\li}{~\\ $\bullet$ }

\newcommand{\di}{\partial}

\newcommand{\smin}{{\setminus}}\newcommand{\sset}{{\subset}}

\newcommand{\Nnd}{Newton-non-degenerate}
\newcommand{\wrt}{with respect to }
\newcommand{\ND}{Newton diagram}

\title{T\MakeLowercase{he} M\MakeLowercase{ilnor fibre signature is not semi-continuous}}
\author{D\MakeLowercase{mitry} K\MakeLowercase{erner} \MakeLowercase{and} A\MakeLowercase{ndr\'as} N\MakeLowercase{\'emethi}}
\address{Department of Mathematics, Ben Gurion University of the Negev, P.O.B. 653, Be'er Sheva 84105, Israel.}
\email{kernerdm@math.bgu.ac.il}
\address{R\'enyi Institute of Mathematics\\Budapest\\Re\'altanoda u. 13\textendash15\\1053\\Hungary}
\email{nemethi@renyi.hu}
\thanks{D.K. was constantly supported by the Skirball postdoctoral fellowship of the Center of Advanced Studies
in Mathematics (Mathematics Department of Ben Gurion University, Israel).}
\thanks{A.N. is partially supported by OTKA and T\'ET Grants of the
Hungarian Academy of Sciences.}

\begin{document}\setcounter{secnumdepth}{6} \setcounter{tocdepth}{6}

\maketitle
\begin{abstract}
Consider the germ of an isolated surface singularity in
$(\mC^3,0)$. The corresponding Milnor fibre possesses the homology
lattice (the integral middle homology with a natural symmetric
intersection form).

An old question of A.Durfee (1978) asks: is the signature of this form
non-increasing under degenerations? The present article answers
negatively: We give examples of Newton non-degenerate families where
the signature increases under degeneration.
\end{abstract}
\section{Introduction}
Let $(X,0)=(f^{-1}(0),0)\sset(\mC^3,0)$ be the germ of a complex,
locally analytic surface with isolated singularity at the origin.
For a general introduction to the `Milnor  package' associated
with $(X,0)$, see \cite{Milnor68}, \cite{Dimca92}, \cite{AGLV} or
\cite{Seade}.
\\
\\
Let $B\sset\mC^3$ be a small ball centered at the origin. Consider a
small deformation: $X_\ep=f^{-1}(\ep)\cap B$, where $\ep\neq0$ is
sufficiently small with respect to  the chosen radius of the ball.
Then $X_\ep$ is smooth and  the diffeomorphism type of the pair
$(B,X_\ep)$ does not depend on $\ep\neq0$. $X_\ep$ is called the
Milnor fibre and it has the homotopy type of a bouquet  $\vee_\mu S^2$
 of middle dimensional spheres, where $\mu$ is the Milnor number,
hence
$H_2(X_\ep,\mZ)\approx\mZ^\mu$. There is a natural symmetric
intersection form on $H_2(X_\ep,\mZ)$. Its extension to
$H_2(X_\ep,\mZ)\otimes\mR$ gives the triple $(\mu_+,\mu_0,\mu_-)$,
the dimensions of maximal vector subspaces of $H_2(X_\ep,\mZ)\otimes\mR$ on which the form is positive definite, zero or negative definite.
\\
\\
Let $(X_t,0)\sset (\mC^3,0)\times(\mC^1_t,0)$ be a holomorphic family of surface germs, such that
$(X_t,0)$ has an isolated singularity for any $t$ in the neighborhood of $0\in\mC^1_t$. An
old question of A.Durfee \cite[Conjecture 5.4]{Durfee}, see also W.Neumann
\cite[Problems session, p.249]{N.T.R.82}, is about the behavior of signature
$sign=\mu_+-\mu_-$ of the corresponding intersection forms under the
degeneration:

{\it Is the signature non-increasing under the degeneration $t\to0$?}
\\
\\
It seems that the question was not addressed previously. The common
expectation was that the semi-continuity as above should hold, at
least under some mild restrictions. And the counterexamples, if
any, must be quite complicated.

Since any non--smooth singularity is realized as the degeneration of
$(x,y,z)\mapsto x^2+y^2+z^2$ whose signature is $-1$,
a positive answer to the above  question would imply that the
signature of any non--smooth isolated hypersurface surface singularity is
negative, which is the (weak) Durfee Conjecture
\cite[Conjecture 5.2]{Durfee}.

The present article shows that any tentative proof of the Durfee Conjecture which is
based on semi--continuity of the  signature in this direct way, will fail.

Note also that the signature can be `almost' determined from the spectrum, which has a rather strong
semi-continuity behavior, cf. \cite{Varchenko83}, \cite{Steenbrink85}. The contribution in the
signature, which is not covered by the spectrum, is the equivariant signature corresponding to the eigenvalue
one of the monodromy, which depends on the structure of the Jordan blocks associated with this
eigenvalue. The point is that exactly these blocks are responsible for the non--semi--continuity
of the signature; and they are the key targets of the present article as well.

It is quite a surprise (at least for us) that the semi-continuity is violated even for families
of {\it \Nnd} surfaces. (In fact, all our examples are degenerations of $T_{p,q,r}$ singularities.)
\\
\bex\label{Example 1}
Consider the family of surfaces $X_t=f^{-1}_t(0)$, where for any $t\in \mC$\\
\parbox{14cm}
{\beq f_t=t xyz+xyz(x+y+z)+x^4y+y^4z+z^4x.
 \eeq Note that $X_t$ is
\Nnd~for any $t$. In fact, $X_{t\not=0}$ is equivalent with the
$T_{13,13,13}$--singularity $\{xyz+x^{13}+y^{13}+z^{13}=0\}$.
Indeed, first by scaling of the
coordinates, one gets that $f_{t\not=0}$ is contact equivalent to $xyz+x^4y+y^4z+z^4x$.
Further, using the transformation
\beq\ber
x\to x+4y^2z^3-6yz^6+4z^9+29y^{11}z^2,\\
y\to x+4x^3z^2-6x^6z+4x^9+29z^{11}x^2,\\
z\to z+4x^2y^3-6xy^6+4y^9+29x^{11}y^2
\eer\eeq
 }
\begin{picture}(0,0)(-20,-20)
\curve(0,0,0,40)\curve(0,0,70,0)\curve(0,0,-35,-35)\put(70,3){$\hat y$}\put(-38,-30){$\hat x$}\put(2,45){$\hat z$}

\curve(-10,30,50,10)\curve(-20,-30,50,10)\curve(-10,30,-20,-30)
\put(-12.5,28){$\bullet$}\put(48.5,8){$\bullet$}\put(-21.5,-31.5){$\bullet$}
\put(-30,37){$\tinyM(0,-1,3)$}\put(42,18){$\tinyM(-1,3,0)$}\put(-26,-40){$\tinyM(3,0,-1)$}

\curve(-10,-10,0,10)\curve(-10,-10,20,0)\curve(0,10,20,0)

\curve(-10,30,-10,-10)\curve(-10,30,0,10)\curve(-20,-30,-10,-10)\curve(-20,-30,20,0)
\curve(50,10,0,10)\curve(50,10,20,0)

\put(-25,-55){\smallT The \ND~ of $\frac{f_0}{xyz}$}
\end{picture}
\eex
~\\
and dividing  by  $1+82(x^8y^2+82y^8z^2+82z^8x^2)$ one gets
$xyz+x^{13}+y^{13}+z^{13}+...$, where the higher order terms consist of monomials of total degree at least 15.

The relevant numerical data for the family is:
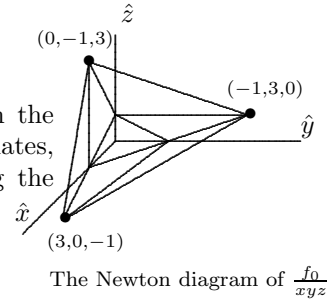
\beq \bM\hline
&\mu&\mu_+&\mu_0&\mu_-&\mu_+-\mu_-\\\hline
X_0&45&5&3&37&-32\\\hline X_{t\neq0}&38&1&1&36&-35\\\hline \eM
\eeq

Once a counterexample is found a natural question is: {\it how significantly the signature
can grow in degenerations}?

Below we consider the families of \Nnd~surface singularities. We
explain the idea behind example \ref{Example 1} and present
additional interesting families. In particular, it appears that
the signature can increase significantly. We give an example when
the change of $\mu_+-\mu_-$ is asymptotically
$\frac{\mu(X_0)}{12}$ (or $\frac{\mu(X_{t\not=0})}{9}$).
\beR If one considers the family $X_t$ as a deformation of the
central fibre and allows the singular point to split into several
ones then an immediate counterexample is obtained as the
suspension of a \Nnd~curve singularity.

Let $(C_0,0)\sset(\mC^2,0)$ be a singularity of the topological
type of $x^3=y^9$, i.e. three smooth branches, pairwise tangent
with order 3.
Consider the family $C_t$ in which the
initial singularity deforms to five singular points, of the types
 $(A_5,A_5,A_1,A_1,A_1)$.
An example of such a family is: $f_t(x,y)=y(y-x^3)(y+(x-t)^3)$.
Here the $A_5$ singularities are at the points $(0,0)$ and
$(t,0)$, while the three nodes $A_1$ are the intersection points
of the branches $\{y-x^3=0\}$ and $\{y+(x-t)^3=0\}$.

Let $X_t=\{z^2+f_t(x,y)=0\}\sset\mC^3$ be the stabilization of
such a family. The relevant data for the signature is:
\beq
\bM\hline &\mu&\mu_+&\mu_0&\mu_-&\mu_+-\mu_-\\\hline
X_0&16&2&0&14&-12\\\hline X_{t\neq0}&13&0&0&13&-13\\\hline \eM
\eeq
(where, for each fixed $t$, the invariants are summed up over
all singular points of $X_t$). Hence $sign(X_0)>sign(X_{t\not=0})$.

  The original question \cite[Conjecture 5.6]{Durfee} and  in
\cite{N.T.R.82} targets non--splitting families, hence from now on
we consider only this case: $Sing(X_t)=\{0\}$ for all $t$.
\eeR~
\\
\\
{\bf Acknowledgements.} The results of this paper were obtained during the conference ``Libgober 60", Jaca,
Spain (July 2009). We appreciate the hard labor of the organizers, excellent conditions and working atmosphere.

For numerical computations we used \cite{Singular}, in particular
the library: "Gau\ss-Manin connection" by M.Schulze. We appreciate the power of the program and the help of the
Singular forum.

Finally, we thank the referee for the careful reading of the manuscript.
\section{Preliminaries}
\subsection{The Milnor number}
For isolated \Nnd~ singularities various topological invariants
are determined by the \ND.  It is defined as follows: first, one takes the
Minkowski sum of the support of $f$ and the first octant $\mR^3_{\ge0}$.
Then the union of all the compact faces of the convex hull of this set  is called the \ND, and it is denoted by
 $\Ga_f$.

The function $f$ is \Nnd~ or non-degenerate with respect to its diagram if for any face $\si\sset\Ga_f$ the restriction $f|_\si$
 is non-degenerate, i.e. the corresponding hypersurface has no singular points inside the maximal torus.

Let $X=f^{-1}(0)\sset(\mC^3,0)$ be a \Nnd~
germ with \ND~ $\Ga_f$. We can assume that $\Ga_f$ is convenient
(by adding to $f$ the monomials $x^N,y^N,z^N$ for $N\gg0$).

Let $\mR^2_{\ge0}\supset\Ga_{xy}=\Ga_f\cap\{z=0\}$, similarly for
$\Ga_{yz}$, $\Ga_{xz}$. Note that all these diagrams are convex
and  convenient. Let
$\mR^1_{\ge0}\supset\Ga_{x}=\Ga_f\cap\{z=0=y\}$, and similarly
introduce  $\Ga_{y}$, $\Ga_{z}$.

The Milnor number is determined by the volumes under $\Ga_f$ as
follows \cite{Kouchnirenko76} \beq\label{Eq. Milnor number
Kouchnirenko} \mu=3! Vol_3-2!Vol_2+Vol_1-1, \eeq where $Vol_3$ is
the volume in $\mR^3_{\ge0}$ under the Newton diagram $\Gamma_f$,
$Vol_2$ is the sum of volumes in $\mR^2_{\ge0}$ under the diagrams
$\Ga_{xy}$, $\Ga_{yz}$, $\Ga_{xz}$, and $Vol_1$ is the sum of
volumes in $\mR^1_{\ge0}$ under the diagrams $\Ga_{x}$, $\Ga_{y}$,
$\Ga_{z}$.  All the volumes are normalized so that the volume of
the unit cube (in any dimension) is 1.
\\
\subsection{$(\mu_+,\mu_0,\mu_-)$}\label{Sec.Preliminaries.Signature via Newton Diagram}
Given a symmetric intersection form let $M_+,M_-,M_0$ be the maximal subspaces of $H_2(X_\ep,\mR)$ on which the symmetric form is positive
definite, zero, negative definite.
 So $\mu_-=dim(M_-)$, and similarly for $\mu_0,\mu_+$.

Recall \cite{Steenbrink76} that the spectrum of $(X,0)$ is a set
of rational numbers in $(-1,2)$, with multiplicities. The spectrum
is symmetric \wrt~$\frac{1}{2}$. The values $Sp\cap (0,1)$
correspond to $M_-$, those from
$\big(Sp\cap(-1,0)\big)\cup\big(Sp\cap (1,2)\big)$ correspond to
$M_+$ and those from $Sp\cap\{0,1\}$ correspond to $M_+\cup M_0$.
In this last case,  the precise distribution into $M_+,M_0$ is
determined by the Jordan block structure of the monodromy (which
is {\em not} codified in the spectrum).
Indeed, the cardinality of $Sp\cap\{0,1\}$ is the dimension of the
generalized eigen-space of the monodromy corresponding to
eigenvalue one. Let $J_n$ denote the number of $n\times n$ Jordan
blocks corresponding to the eigenvalue one. By the monodromy theorem $n$ is always $\leq
2$, see e.g. \cite{St} and the references therein.
Then $J_1$ contributes 1 to $\mu_0$ (and provides either the
spectral number 0 {\em or} the spectral number 1; by the symmetry
of the spectrum such blocks appear in `pair', one of them having
spectral number 0, the other 1), while $J_2$ contributes 1 to
$\mu_0$ and 1 to $\mu_+$ (and provides two spectral numbers,
namely  0 {\em and} 1).

The cardinalities of the spectral numbers in different intervals (as
above) are computed from the \ND~as follows \cite{Saito88}:

\li $\sharp\big(Sp\cap(-1,0)\big)$ equals the number of
$\mZ^3_{>0}$ points strictly below $\Ga_f$.

\li $\sharp \big(Sp\cap\{0\}\big) $  equals the number of
$\mZ^3_{>0}$ points on $\Ga_f$; hence:

 \li $\sharp\big(Sp\cap(-1,0]\big)$ equals the
number of $\mZ^3_{>0}$ points not above $\Ga_f$. Hence,
$\mu_0+\mu_+$ is twice this number.

\vspace{2mm}

\noindent Regarding the Jordan block structure on the eigenvalue 1
part of the monodromy one has (\cite{St}, p. 227):

\li $J_2$ is the number of  points of $\mZ^3_{>0}\cap\Ga_f$ which
do not  lie  in the interior of a two-dimensional (compact) face
of $\Ga_f$.

\bex Let $X=f^{-1}(0)\sset\mC^3$ for $f=xyz+x^p+y^q+z^r$. There is
only one point of $\mZ^3_{>0}$ which does not lie  above $\Ga_f$, namely
$(1,1,1)$; hence $\mu_++\mu_0=2$. This point does not lie in the
interior of any face, hence it contributes to $J_2$ and
$\mu_0=1=\mu_+$.

Applying Kouchnirenko's formula (or in any other way) we get: $\mu=p+q+r-1$, $\mu_-=p+q+r-3$
and the signature $\mu_+-\mu_-=4-p-q-r$.
\eex
\subsection{Numerical computations}
The general procedures to compute the spectrum and the monodromy are realized in \cite{Singular} (though they
are often time consuming).
For, example the code for the case above is:
\\
\\LIB "sing.lib";\hspace{6cm}\%\%\%\%  Loading the basic library of singularities.
\\LIB "gmssing.lib";\hspace{5.4cm}\%\%\%\%  Loading the library of Gau{\ss}-Manin connection.
\\ring s = 0,(x,y,z),ds;\hspace{5.2cm}\%\%\%\%  Local ring in 3 variables is defined
\\int p,q,r=4,5,6;\hspace{6cm}\%\%\%\%  Integer constants
\\poly $f=x*y*z+x^p+y^q+z^r$;
\\milnor(f);
\\spectrum(f);
\\monodromy(f);
\section{The idea behind the examples}
As we consider the families of \Nnd~singularities, we prescribe the
degeneration by erasing some vertices.
For simplicity, we will only describe a `primitive' degeneration:
just one vertex is erased. By assumption, the vertex does not
belong to the interior of a face or an edge. Hence, it is enough
to consider just the relevant part of the \ND.
\\
\\
\parbox{12cm}
{Let $a\in\mZ^3_{>0}$ be the apex of a (strictly convex, solid)
cone $C$.
Erase the apex, get a new smaller body:
$C_{new}=Conv\big((C\setminus a)\cap\mZ^3_{\ge0}\big)$. Let $\di
C$, respectively  $\di C_{new}$, denote the boundary of $C$,
respectively $C_{new}$. An integral point of $\di C_{new}\smin
\,\di C$ is called {\it new}.

We assume that there is no new point
which lies in the interior of a 2--face of $C_{new}$.
Consider the boundary of $\di C_{new}\setminus \di C$, it is a
closed 1--dimensional piecewise linear loop. The interior points
of its edges are called {\it inner}. Its  vertices  are called
{\it outer}. Finally, let $V$ denote the volume of $C\setminus
C_{new}$.

}
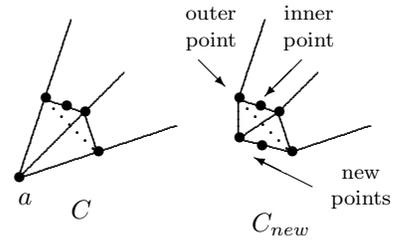
\begin{picture}(0,0)(-20,20)
\curve(0,0,10,30,20,60)\curve(0,0,20,20,40,40)\curve(0,0,30,10,60,20)
\curve(10,30,25,25)\curve(30,10,25,25)\multiput(9,29)(3,-3){7}{.}
\put(16,25){$\bullet$}
\put(-2,-2){$\bullet$}\put(8,28){$\bullet$}\put(28,8){$\bullet$}\put(23,23){$\bullet$}

\put(0,-10){$a$}
\put(20,-15){$C$}
\end{picture}
\begin{picture}(0,0)(-90,20)
\curve(10,30,20,60)\curve(25,25,40,40)\curve(30,10,60,20)
\curve(10,30,25,25)\curve(30,10,25,25)\multiput(9,29)(3,-3){7}{.}
\curve(10,15,10,30)\curve(10,15,30,10)\curve(10,15,25,25)

\put(8,13){$\bullet$}\put(16.5,10){$\bullet$}\put(8,28){$\bullet$}\put(28,8){$\bullet$}\put(23,23){$\bullet$}
\put(16,25){$\bullet$}
\put(15,-20){$C_{new}$}

\put(38,-3){\vector(-2,1){22}}\put(40,-5){\smallT \bet new\\points\eet}

\put(-5,45){\vector(1,-1){10}}\put(-15,55){\smallT \bet outer\\point\eet}
\put(34,45){\vector(-1,-1){14}}\put(22,55){\smallT \bet inner\\point\eet}
\end{picture}
\\
\\
\\
 Using these notations and assumptions, the  discussion of
 (\ref{Sec.Preliminaries.Signature via Newton Diagram}) implies
 the following:
\bprop
Let $N_{inner}$, $N_{outer}$, $N_{new}$ be the number of inner,
outer and new points. Then: \beq\ber
\mu(X_0)-\mu(X_{t\neq0})=6V\\
\mu_0(X_0)-\mu_0(X_{t\neq0})=-1-N_{inner}+N_{new}\\
\mu_+(X_0)-\mu_+(X_{t\neq0})=1+N_{inner}+N_{new}\\
\mu_-(X_0)-\mu_-(X_{t\neq0})=6V-2N_{new}\\
sign(X_0)-sign(X_{t\neq0})=-6V+N_{inner}+3N_{new}+1. \eer\eeq
\eprop

Next, we estimate $V$. Take an arbitrary elementary integral
triangulation of the 2-faces of $\di C_{new}\smin\,\di C$, i.e. the
vertices of each triangle belong to $\mZ^3_{\ge0}$
 and there are no integral points in the interior
 of the triangle or of its edges.
Then $C\smin C_{new}$  is naturally subdivided into the pyramids
with the triangles as bases and the apex $a\in C$. Each pyramid is
elementary, i.e. its only integral points are its corners.
\bel
Set $\delta:=6V - N_{outer}-N_{inner}-2N_{new}+2$. Then $\delta
\ge 0$. The equality occurs iff the volume of each of the pyramids
above is $\frac{1}{6}$.
\eel
\bpr By the construction above we get
that $6V\ge$the number of pyramids=the number of triangles in $\di
C_{new}\smin\, \di C$. The later quantity is combinatorial, does
not depend on a particular triangulation, and can be  computed
explicitly.
\epr

From the previous proposition we get
$sign(X_0)-sign(X_{t\neq0})=3-\delta-N_{outer}+N_{new}$.
Therefore, unsing this lemma in order to get an interesting example (when this
expression is `large', or at least it is positive),  $N_{outer}$
should be small, $N_{new}$ should be large and the pyramids above
should have the minimal volume $\frac{1}{6}$ (or as small as
possible).
\section{A degeneration with big increase of signature}
A simple generalization of the example \ref{Example 1} is the
following.
\bprop The invariants for the family of surfaces
$(X_t,0)=f^{-1}_t(0)\sset(\mC^3,0)$ with
\\
\parbox{13cm}
{\beq f_t=txyz+x^{3k+3}y+x^{k+1}yz+z^2x+y^2z
\eeq are:
\beq
\bM\hline &\mu&\mu_+&\mu_0&\mu_-&\mu_+-\mu_-\\\hline
X_0&12k+11&2k+1&1&10k+9&-8k-8\\\hline
X_{t\neq0}&9k+11&1&1&9k+9&-9k-8\\\hline \eM
\eeq So the signature
grows by $k$, therefore asymptotically by $\frac{\mu(X_0)}{12}$ (or
$\frac{\mu(X_{t\not=0})}{9}$) when $k\mapsto\infty$. }
\begin{picture}(0,0)(-80,-20)
\curve(0,0,0,40)\curve(0,0,30,0)\curve(0,0,-60,-40)\put(25,-10){$\hat y$}\put(-67,-37){$\hat x$}\put(-2,45){$\hat z$}

\curve(-5,30,20,5)\curve(-5,30,-24,-16)\curve(-5,30,-50,-40)
\curve(20,5,-24,-16)\curve(-50,-40,-24,-16)\curve(-50,-40,20,5)

\put(-52,-42){$\bullet$}\put(-7,28){$\bullet$}\put(18,3){$\bullet$}\put(-26,-18){$\bullet$}
\put(-65,-47){$\tinyA(3k+2,0,-1)$}\put(-30,35){$\tinyA(0,-1,1)$}\put(18,10){$\tinyA(-1,1,0)$}

\put(-65,-65){\smallT The \ND~ of $\frac{f_t}{xyz}$}
\end{picture}
\eprop
Notice that $X_{t\not=0}$ is equivalent to the singularity $T_{6k+5,3k+4,3}$.
This can be proved in many different ways, and, in fact, all these methods can be used in the
 similar statement of Example \ref{Example 1} as well.
 By the first method, one computes the resolution graph of the germ, e.g. by Oka's algorithm \cite{Oka87}
 (cf. also \cite{BN}), and one gets a cyclic graph which characterizes the hypersurface $T_{p,q,r}$ germs.
 Or, one can use the results \cite{BN} on `equivalent Newton boundaries' which identifies
 the triple $(6k+5,3k+4,3)$ from the diagram.  It is clear that $f_{t\not=0}$ is equivalent to
 $xyz+x^{3k+3}y+z^2x+y^2z$. Then consider the three 2--faces of the the diagram and notice that they
 intersect the three axes at $(6k+5,0,0)$, $(0,3k+4,0)$ and $(0,0,3)$ respectively. Then use
 \cite[\S 3]{BN}. Finally, one can get the result by change of variables too.

\bpr
For the computation of $\mu_0$, $\mu_+$ we consider the polynomial $\frac{f}{xyz}$, then we look for the point of $\mZ^2_{\ge0}$
 lying under/not above $\Ga_{\frac{f}{xyz}}$.
One sees that any such point is of the form $(x,0,0)$. In addition
no integral points lie in the interior of the faces or edges of
$\Ga_f$ (or $\Ga_{\frac{f}{xyz}}$). Hence: \li $t\neq0$. There is
only one point: $(0,0,0)\in\Ga_{\frac{f}{xyz}}$ it creates $J_2$,
so contributes 1 to $\mu_0$ and 1 to $\mu_+$. \li $t=0$. The
points $(j,0,0)$ for $0\le j<k$ lie under $\Ga_{\frac{f}{xyz}}$
they contribute $2k$ to $\mu_+$. The point
$(k,0,0)\in\Ga_{\frac{f}{xyz}}$  contributes 1 to $\mu_0$ and 1 to
$\mu_+$.

The Milnor number can be computed by Kouchnirenko formula or in any other way, then one obtains $\mu_-$ too.
\epr

\end{document}